\documentclass[11pt]{amsart}
\usepackage{amsmath,amssymb}
\usepackage{amscd,eucal,amsthm}
\newtheorem{Theorem}{Theorem}
\newtheorem{Corollary}{Corollary}
\newtheorem{Lemma}{Lemma}
\newtheorem{Fact}{Fact}
\newtheorem{Proposition}{Proposition}

\theoremstyle{remark}
\newtheorem*{Remark}{Remark}

\renewcommand{\to}[1][]{\xrightarrow{#1}}

\newcommand{\C}{{\mathbb{C}}}

\newcommand{\Z}{{\mathbb{Z}}}
\renewcommand{\gg}{{\mathfrak{g}}}
\renewcommand{\ll}{{\mathfrak{l}}}
\newcommand{\nn}{{\mathfrak{n}}}

\newcommand{\z}{{\mathfrak{z}}}

\renewcommand{\ss}{{\mathfrak{s}}}
\newcommand{\bb}{{\mathfrak{b}}}
\newcommand{\hh}{{\mathfrak{h}}}
\renewcommand{\l}{{\lambda}}

\newcommand{\rr}{\mathbb{R}}

\newcommand{\cc}{\mathbb{C}}

\newcommand{\PP}{\mathbb{P}}
\newcommand{\A}{\mathcal{A}} 
\DeclareMathOperator{\can}{can} \DeclareMathOperator{\Op}{Op}
 \DeclareMathOperator{\re}{Re}
 \DeclareMathOperator{\gr}{gr}
\DeclareMathOperator{\ad}{ad} \DeclareMathOperator{\Ad}{Ad}

 \DeclareMathOperator{\Arg}{Arg}
 
\DeclareMathOperator{\tr}{tr} \DeclareMathOperator{\End}{End}
 
 \DeclareMathOperator{\Spec}{Spec}
 \DeclareMathOperator{\rk}{rk}
 
\renewcommand{\phi}{\varphi}

\newcommand{\ka}{\kappa}
\newcommand{\la}{\lambda}

 \makeatletter
\def\@mult#1{\raise #1\rlap{$\cdot$}\lower #1\rlap{$\cdot$}\cdot}
\def\did{\mathrel{\@mult{3pt}}}
\makeatother
\def\openrow#1#2#3{\setbox0=\vbox{\hbox
    {\vrule height#2 width#3\kern#2\vrule height#2 width0pt}\hrule height#3}
    \hbox{\leaders\copy0\hskip#1\wd0\vrule width#3}}
\def\row#1#2#3{\vbox{\hrule height#3\openrow{#1}{#2}{#3}}}
\def\Yr#1{\row{#1}{1.5ex}{.1ex}}

\def\DY#1\endDY{\baselineskip=1ex\lineskip=0pt\lineskiplimit=0pt{\vcenter
    {\Yr#1}}}
\def\openclm#1#2#3{\setbox0=\vbox{\hrule height#3\hbox
    {\vrule width0pt\kern#2\vrule width#3 height#2}}\vtop
    {\leaders\copy0\vskip#1\ht0\hrule height#3}}
\def\clm#1#2#3{\hbox{\vrule width#3\openclm{#1}{#2}{#3}}}
\def\Yc#1{\clm{#1}{1.5ex}{.1ex}}

\def\CDY#1\endCDY{{\vcenter{\hbox{\Yc#1}}}}

\newcommand{\nc}{\newcommand}
\nc{\g}{{\mathfrak g}}
\nc{\LG}{{}^L\neg G}
\nc{\pone}{{\mathbb P}^1}
\nc{\wt}{\widetilde}
\nc{\wh}{\widehat}
\nc{\ghat}{\wh{\gg}}
\nc{\mc}{\mathcal}
\nc{\su}{{\mathfrak s}{\mathfrak l}}
\nc{\ppart}{(\!(t)\!)}
\nc{\on}{\operatorname}
\nc{\sw}{{\mf s}{\mf l}}
\nc{\mf}{\mathfrak}
\nc{\ol}{\overline}
\nc{\Gr}{\on{Gr}}
\nc{\bi}{\bibitem}

\begin{document}

\author[Boris Feigin]{Boris Feigin$^1$}

\thanks{$^1$ Supported by the grants RFBR 05-01-01007, RFBR
05-01-02934 and NSh-6358.2006.2.}

\address{Landau Institute for Theoretical Physics, Kosygina St 2,
Moscow 117940, Russia}

\author[Edward Frenkel]{Edward Frenkel$^2$}

\thanks{$^2$ Supported by DARPA and
AFOSR through the grant FA9550-07-1-0543.}

\address{Department of Mathematics, University of California,
  Berkeley, CA 94720, USA}

\author[Leonid Rybnikov]{Leonid Rybnikov$^3$}

\thanks{$^3$ Supported by
the grants RFBR 05 01 00988-a, RFBR 05-01-02805-CNRSL-a, and
Deligne 2004 Balzan prize in mathematics. The work was finished
during L.R.'s stay at the Institute for Advanced Study supported
by the NSF grant DMS-0635607.}

\address{Poncelet laboratory (IUM and CNRS) and Institute for
Theoretical and Experimental Physics, Moscow, Russia}

\title{Opers with irregular singularity and spectra of the shift of argument
subalgebra}

\date{December 2007}

\begin{abstract}
The universal enveloping algebra of any simple Lie algebra
$\g$ contains a family of commutative subalgebras, called the quantum
shift of argument subalgebras \cite{Ryb1,FFTL}. We prove that
generically their action on finite-dimensional modules is
diagonalizable and their joint spectra are in bijection with the set
of monodromy-free $^L G$-opers on $\pone$ with regular singularity at
one point and irregular singularity of order two at another point. We
also prove a multi-point generalization of this result, describing the
spectra of commuting Hamiltonians in Gaudin models with irregular
singulairity. In addition, we show that the quantum shift of argument
subalgebra corresponding to a regular nilpotent element of $\g$ has a
cyclic vector in any irreducible finite-dimensional $\g$-module. As a
byproduct, we obtain the structure of a Gorenstein ring on any
such module. This fact may have geometric significance related to the
intersection cohomology of Schubert varieties in the affine
Grassmannian.
\end{abstract}

\maketitle

\section{Introduction}

Let $\g$ be a simple Lie algebra over $\cc$. The symmetric algebra
$S(\g)$ carries a natural Poisson structure. A Poisson commutative
subalgebra $\overline{\mc A}_\mu$ of $S(\g) = \cc[\g^*]$, called
the {\em classical shift of argument subalgebra}, was defined in
\cite{MF} (see also \cite{Man}). It is generated by the
derivatives of all orders in the direction of $\mu\in\gg^*$ of all
elements of the algebra of $\g$-invariants in $\cc[\g^*]$.
Recently, this algebra was quantized in \cite{Ryb1,FFTL}. More
precisely, a commutative subalgebra ${\mc A}_\mu$ of the universal
enveloping algebra $U(\gg)$ of $\gg$ was constructed and it was
proved that for all regular $\mu\in\gg^*$ the associated graded of
${\mc A}_\mu$, with respect to the order filtration on $U(\g)$ is
$\overline{\mc A}_\mu$. The algebra ${\mc A}_\mu$ is called the
{\em quantum shift of argument subalgebra} of $U(\g)$. This is a
free polynomial algebra in $\frac{1}{2}(\dim\gg+\rk\gg)$
generators for any regular $\mu \in \gg^*$ (see \cite{Ryb1,FFTL}).

For any regular $\mu$ the algebra ${\mc A}_\mu$ contains the
centralizer of $\mu$ in $\gg$. In particular, if $\mu$ is regular
semi-simple, then ${\mc A}_\mu$ contains a Cartan subalgebra of
$\gg$. Hence it acts on weight subspaces of $\g$-modules. We note
that, as shown in \cite{Ryb1}, a certain limit of ${\mc A}_\mu$ in
the case when $\g=\sw_n$ may be identified with the
Gelfand--Zetlin algebra. Hence the algebra ${\mc A}_\mu$ may be
thought of as a generalization of the Gelfand--Zetlin algebra to
an arbitrary simple Lie algebra.  It is an interesting question to
describe the joint (generalized) eigenvalues of ${\mc A}_\mu$ on
$\g$-modules.

The first steps towards answering this question were taken in
\cite{FFTL}, where it was shown that ${\mc A}_\mu$ is isomorphic to
the algebra of functions on a certain space of $^L G$-{\em opers} on
the projective line $\PP^1$. More precisely, this is the space of $^L
G$-opers on $\PP^1$ with regular singularity at the point $0$ and
irregular singularity of order $2$ at the point $\infty$, with a fixed
$2$-residue determined by $\mu$ (see \cite{FFTL} for details). Here
$^L G$ is the {\em Langlands dual group} of $G$ ($^L G$ is taken to be
of adjoint type), and $^L G$-opers are connections on a principal $^L
G$-bundle over $\pone$ satisfying a certain transversality condition,
as defined in \cite{BD}. The appearance of the Langlands dual group is
not accidental, but is closely related to the geometric Langlands
correspondence, through a description of the center of the completed
enveloping algebra of the affine Kac--Moody algebra $\ghat$ at the
critical level in terms of $^L G$-opers on the punctured disc
\cite{FF,Fr2,F:book}.

Thus, we obtain that the spectra of ${\mc A}_\mu$ on a $\g$-module
$M$ are encoded by $^L G$-opers on $\pone$ satisfying the above
properties. Furthermore, in \cite{FFTL} it was shown that if
$M=V_\la$, the irreducible finite-dimensional $\g$-module with
dominant integral highest weight $\la$, then these $^L G$-opers
satisfy two additional properties: they have a fixed residue at
the point $0$ (where the oper has regular singularity), determined
by $\la$ (we denote the space of such opers by $\Op_{^L
G}(\PP^1)_{\pi(-\mu)}^\l$), and they have {\em trivial monodromy}.
It was conjectured in \cite{FFTL} that in fact there is a
bijection between the spectra of ${\mc A}_\mu$ on $V_\la$ and the
set of monodromy-free opers from $\Op_{^L
G}(\PP^1)_{\pi(-\mu)}^\l$.

In this paper we prove this conjecture. Furthermore, we prove the
following statement (see Corollary~\ref{simple-spectrum}):

\medskip

\noindent {\bf Theorem A.} {\em For generic regular semi-simple
$\mu\in\gg^*$ and any dominant integral $\l$, the quantum shift of
argument subalgebra $\A_\mu \subset U(\gg)$ is diagonalizable and
has simple spectrum on the irreducible $\g$-module $V_\l$.
Moreover, its joint eigenvalues (and hence eigenvectors, up to a
scalar) are in one-to-one correspondence with monodromy-free opers
from $\Op_{^L G}(\PP^1)_{\pi(-\mu)}^{\l}$.}

\medskip

In this paper we do not address the question of constructing the
eigenvectors of ${\mc A}_\mu$. Conjecturally, for generic $\mu$, they
may be constructed by the Bethe Ansatz method described in
\cite{FFTL}, but we do not attempt to prove this conjecture here, nor
do we use Bethe Ansatz in the proof of Theorem A. Instead, we rely on
the isomorphism between the algebra ${\mc A}_\mu$ and the algebra of
functions on opers and the study of opers with irregular singularity
and trivial monodromy.

The crucial step in our proof is the analysis of the action of ${\mc
A}_\mu$ in the case when $\mu = f$, a regular {\em nilpotent} element
of $\gg^* \simeq \gg$. The algebra ${\mc A}_f$ contains the
centralizer ${\mf a}_f$ of $f$ in $\gg$. Elements of ${\mf a}_f$ and
more general elements of ${\mc A}_f$ are not diagonalizable operators
on $V_\la$, but nilpotent operators. Hence it is natural to view them
as ``creation operators'' and ask whether they generate $V_\la$ from
its highest weight vector. Theorem~\ref{regular_sequence} implies that
the answer to this question is affirmative (see
Corollary~\ref{monodr_cond} for more details).

\medskip

\noindent {\bf Theorem B.} {\em The $\gg$-module $V_\l$ is cyclic as
an $\A_f$-module. The annihilator of $V_\l$ in $\A_{f}$ is generated
by the no-monodromy conditions on opers from $\Op_{^L
G}(\PP^1)_{0}^\l$.}

\medskip

This provides a natural structure of a Gorenstein ring on any
finite-dimensional irreducible $\gg$-module $V_\l$. We hope that this
fact has a geometric interpretation. Namely, due to the results of
\cite{Gi,MV}, the space $V_\l$ may be naturally identified with the
global cohomology of the irreducible perverse sheaf on the Schubert
variey $\ol\Gr_\l$ in the affine Grassmannian of the Langlands dual
group $^L G$. On the other hand, the universal enveloping algebra
$U({\mf a}_f)$ of the centralizer ${\mf a}_f$ of the principal
nilpotent element $f$ is identified with the cohomology ring of the
affine Grassmannian. Hence it acts on the cohomology of our perverse
sheaf. This action is precisely the action of $U({\mf a}_f)$ on $V_\l$
\cite{Gi}. But $U({\mf a}_f)$ is a subalgebra of the commutative
algebra $\A_f$. This leads us to a natural question: what is the
geometric meaning of $\A_f$? Perhaps, if one could answer this
question, one could derive the cyclicity of $V_\l$ as a $\A_f$-module
by geometric means.

\medskip

The shift of argument subalgebra ${\mc A}_\mu$ has a multi-point
generalization, denoted by ${\mc A}_\mu(z_1,\ldots,z_N)$, where
$z_1,\ldots,z_N$ are distinct points on $\PP^1 \backslash
\infty$. This is a commutative subalgebra of $U(\gg)^{\otimes N}$,
which consists of the Hamiltonians of the {\em Gaudin model} with
irregular singularity \cite{Ryb1,FFTL}. We show that this algebra is
isomorphic to the algebra of functions on the space of $^L G$-opers on
$\PP^1$ with regular singularities at $z_1,\ldots,z_N$ and irregular
singularity of order $2$ at $\infty$, with the $2$-residue determined
by $\mu$. We also prove a conjecture of \cite{FFTL} that the set of
joint eigenvalues of this algebra on the tensor product $V_{\l_1}
\otimes \ldots \otimes V_{\l_N}$ of irreducible finite-dimensional
$\gg$-modules is in bijection with the set of opers of this kind with
fixed residues determined by the highest weights $\l_1,\ldots,\l_N$
and the no-monodromy condition (see
Corollary~\ref{simple-spectrum_z}). This is a multi-point
generalization of Theorem A. We also prove a multi-point analogue of
Theorem B (see Corollary~\ref{regular-elements-z}).

\medskip

We note that for the ordinary Gaudin model (with regular singularity)
a description of the spectrum of the Hamiltonians in terms of
an appropriate set of monodromy-free opers analogous to Theorem A was
conjectured in \cite{F:faro}, Conjecture 1 (it was proved in
\cite{F:faro}, Theorem 2.7,(3) that the spectrum does embed into this
set of monodromy-free opers). In the course of writing this paper we
learned that a variant of this conjecture was proved in \cite{MTV} in
the case when $\gg=\gg\ll_M$, by a detailed analysis of intersections
of Schubert varieties in the Grassmannians. Theorem A in the case of
$\gg=\gg\ll_N$ should be related to this statement via the duality of
\cite{TL}.

\medskip

Finally, we expect that the results of this paper may be generalized
to the affine Kac--Moody algebras. As explained in \cite{FF:sol}, the
affine analogue of the shift of argument subalgebra corresponding to
regular semi-simple $\mu$ is the algebra of quantum integrals of
motion of the AKNS hierarchy of soliton equations. Such an algebra may
also be defined for a regular nilpotent $\mu$. In this case, its
action is not diagonalizable, but it gives rise to a commutative
algebra of creation operators. We expect that these operators generate
highest weight modules over affine Kac--Moody algebras; for example,
the irreducible integrable representations. In the latter case, we
expect that the generators of the corresponding annihilating ideal are
given by the no-monodromy conditions on the corresponding affine
opers, by analogy with Theorem B in the finite-dimensional case. We
plan to discuss this in more detail in our next paper.

\medskip

The present paper is organized as follows. In
section~\ref{sect-prelim} we collect basic facts on Gaudin models and
opers. In section~\ref{sect-results} we formulate and discuss the main
result of the paper. The detailed proof is given in the last
section~\ref{sect-Sibuya}.

\medskip

\noindent {\bf Acknowledgements.} We thank A.~Glutsyuk, Yu.~Ilyashenko,
and V. Toledano Laredo for useful discussions.

\section{Preliminaries}\label{sect-prelim}

\subsection{Gaudin algebras} Gaudin model was introduced in \cite{G1}
as a spin model related to the Lie algebra $\su_2$, and
generalized to the case of an arbitrary semisimple Lie algebra in
\cite{G}, Section 13.2.2. For any $x\in\gg$, set
$$x^{(i)}=1\otimes\dots\otimes 1\otimes x\otimes
1\otimes\dots\otimes 1\in U(\gg)^{\otimes N}$$ ($x$ at the $i$th
place). Let $\{x_a\},\ a=1,\dots,\dim\gg$, be an orthonormal basis
of $\gg$ with respect to Killing form, and let $z_1,\dots,z_N$ be
pairwise distinct complex numbers. The Hamiltonians of Gaudin
model are the following mutually commuting elements of
$U(\gg)^{\otimes N}$:
\begin{equation}\label{quadratic}
H_i=\sum_{k\neq i}\sum_{a=1}^{\dim\gg}
\frac{x_a^{(i)}x_a^{(k)}}{z_i-z_k}.
\end{equation}

In \cite{FFR}, a large commutative subalgebra $\A(z_1,\dots,z_N)$
containing the $H_i$'s was constructed with the help of the affine
Kac--Moody algebra $\ghat$, which is the universal central extension
of $\g\ppart$.  Namely, according to \cite{FF,Fr2}, the completed
enveloping algebra $\wt{U}_{\ka_c}(\wh{\gg})$ of $\ghat$ at the
critical level (in the notation of \cite{FFTL}) contains a large
center $Z(\wh{\gg})$. Set $\wh{\gg}_+=\gg[[t]]\subset\wh{\gg}$ and
$\wh{\gg}_-=t^{-1}\gg[t^{-1}]\subset\wh{\gg}$. The natural
homomorphism $Z(\wh{\gg})\to
(\wt{U}_{\ka_c}(\wh{\gg})/\wt{U}_{\ka_c}(\wh{\gg}) \cdot
\wh{\gg}_+)^{\wh{\gg}_+}$ is surjective \cite{FF,Fr2}. Every
element of the latter quotient has a unique representative in
$U(\wh{\gg}_-)$. Thus we obtain the following natural embedding
$$
(\wt{U}_{\ka_c}(\wh{\gg})/\wt{U}_{\ka_c}(\wh{\gg}) \cdot
\wh{\gg}_+)^{\wh{\gg}_+} \hookrightarrow U(\wh{\gg}_-).
$$
Let ${\mathfrak z}(\ghat) \subset U(\wh{\gg}_-)$ be the image of this
embedding. The commutative subalgebra $\A(z_1,\dots,z_N) \subset
U(\gg)^{\otimes N}$ is then the image of ${\mathfrak z}(\ghat) \subset
U(\wh\gg_-)$ under the homomorphism
\begin{equation}    \label{evaluation}
U(\wh\gg_-)\to U(\gg)^{\otimes N}
\end{equation}
of evaluation at the points $z_1,\dots,z_N$ (see \cite{FFR}).

This construction is generalized as follows (see \cite{FFTL,Ryb1}.
For different approach in the case $\gg=\gg\ll_N$ see \cite{ChT}).
One constructs a family of homomorphisms
\begin{equation}    \label{eval}
U(\wh{\gg}_-)\to U(\gg)^{\otimes N}\otimes S(\gg)
\end{equation}
generalizing the evaluation homomorphisms \eqref{evaluation} in
the sense that composing it with the natural augmentation $S(\g)
\to \C$ we obtain \eqref{evaluation}. Now, for any collection
$z_1,\dots,z_N$, the image of ${\mathfrak z}(\ghat)$ under this
homomorphism is a certain commutative subalgebra $U(\gg)^{\otimes
N}\otimes S(\gg)$. Evaluating at any point $\mu\in\gg^*=\Spec
S(\gg)$, we obtain a commutative subalgebra
$\A_{\mu}(z_1,\dots,z_N)\subset U(\gg)^{\otimes N}$ depending on
$z_1,\dots,z_N$ and $\mu\in\gg^*$. These subalgebras contain the
following ``inhomogeneous'' Gaudin hamiltonians:
$$
H_i=\sum\limits_{k\neq i}\sum\limits_{a=1}^{\dim\gg}
\frac{x_a^{(i)}x_a^{(k)}}{z_i-z_k}+
\sum\limits_{a=1}^{\dim\gg}\mu(x_a)x_a^{(i)}.
$$
In particular, $\A(z_1,\dots,z_N) = \A_{0}(z_1,\dots,z_N)$
corresponding to $\mu=0$.

In particular, for $N=1$, we obtain a family of commutative
subalgebras $\A_{\mu}(z_1) \subset U(\gg)$ which does not depend
on $z_1$.  We will set $z_1=0$ in this case and denote this
algebra simply by ${\mc A}_\mu$. It is proved in \cite{FFTL,Ryb1}
that the associated graded algebra of $\A_\mu$ (with respect to
the PBW filtration) for regular $\mu$ is the (classical) {\em
shift of argument subalgebra} $\overline\A_\mu\subset S(\gg)$.
This Poisson commutative subalgebra was first constructed by
Mishchenko and Fomenko in \cite{MF}, in the following way. Let
\begin{equation}    \label{ZS}
S(\gg)^{\gg}=\cc[P_1,\dots,P_\ell],
\end{equation}
where $\ell = \on{rank}(\gg)$ be the center of $S(\gg)$ with respect
to the Poisson bracket, and the $P_i$ are chosen so that they are
homogeneous with respect to the natural grading on $S(\gg)$. Let
$\mu\in\gg^*$ be a regular semisimple element. Then the subalgebra
$\overline\A_\mu\subset S(\gg)$ is generated by the elements
$\partial_{\mu}^nP_k$, where $k=1,\dots,\ell,\ n=0,\dots,\deg P_k-1$,
(or, equivalently, generated by central elements of
$S(\gg)=\cc[\gg^*]$ shifted by $t\mu$ for all $t\in\cc$). These
elements are algebraically independent (for a unform proof, see
\cite{FFTL}, Theorem 3.11). Hence the subalgebra
$\overline\A_\mu\subset S(\gg)$ is a free polynomial algebra in
$\frac{1}{2}(\dim\gg+\rk\gg)$ generators (and therefore has maximal
possible transcendence degree). Since $\on{gr} {\mc A}_\mu = \ol{\mc
A}_\mu$, we obtain the following

\begin{Lemma}    \label{free pol}
The algebra ${\mc A}_\mu$ is a free polynomial algebra in
$\frac{1}{2}(\dim\gg+\rk\gg)$ generators.
\end{Lemma}

Let
$$
\gg = \nn \oplus \hh \oplus \nn_-
$$
be a Cartan decomposition of $\gg$. Denote by $\gg_\rr$ be the compact
real form of $\gg$, and let $\hh_\rr=\gg_\rr\cap\hh$. Any irreducible
finite-dimensional $\g$-module has a $\gg_\rr$-invariant Hermitian
form. We shall use the following

\begin{Lemma}\label{hermitian}
For $z_1,\dots,z_N\in\rr$ and $\mu\in i\hh^*_\rr \subset \hh^*$,
the algebra $\A_{\mu}(z_1,\dots,z_N)\subset U(\gg)^{\otimes N}$
acts by Hermitian operators on any irreducible finite-dimensional
$\g$-module.
\end{Lemma}

\begin{proof}
This follows from the fact that the center at the critical level
and the homomorphism \eqref{eval} are defined over $\rr$, and
hence $\A_{\mu}(z_1,\dots,z_N)\subset U(\gg_\rr)^{\otimes N}$.
\end{proof}

We also need the following facts on the limit points of the family
$\A_{\mu}(z_1,\dots,z_N)$. Namely,

\begin{Proposition}\label{predel_z}\cite{Ryb1}
$\lim\limits_{s\to\infty}\A_{\mu}(sz_1,\dots,sz_N)=\lim\limits_{s\to\infty}\A_{s\mu}(z_1,\dots,z_N)=
\A_{\mu}^{(1)}\otimes\dots\otimes\A_{\mu}^{(N)}\subset
U(\gg)^{\otimes N}$ for regular semisimple $\mu\in\gg^*$.
\end{Proposition}

We have $\A_\xi=\A_{t\xi}$ for any $t\in\cc^\times,
\xi\in\gg^*$. Thus, the subalgebras $\A_\xi\subset U(\gg)$ and
$\A_{t\Ad(g)\xi}\subset U(\gg)$ are conjugate for any $g\in G,
t\in\cc^*, \xi\in\gg$.  Let $\Pi$ be the set of simple roots of
$\gg$. We will choose generators $\{ e_{-\alpha} \}_{\alpha \in \Pi}$
of the lower nilpotent subalgebra $\nn_-$. Let
$$
f=\sum_{\alpha\in\Pi}e_{-\alpha}\in\gg
$$
be a principal nilpotent element. From now on we will identify $\gg$
and $\gg^*$ using a non-degenerate invariant inner product.

\begin{Lemma}\label{predel_mu}
For any regular $\mu\in\gg$, the closure of the family
$\A_{t\Ad(g)\xi}\subset U(\gg)$ ($g\in G, t\in\cc^*, \xi\in\gg$)
contains the subalgebra $\A_f\subset U(\gg)$.
\end{Lemma}

\begin{proof}
Let $\{ e,\ h,\ f \}$ be a principal $\sw_2$-triple in $\gg$, where $h
\in \hh$. Consider Kostant's slice of regular elements
$$
\gg_{\can}:=f+\z_\gg(e).
$$
For any regular $\mu\in\gg$, the $\Ad G$-orbit of $\mu$ intersects
$\gg_{\can}$ (this is a classical result due to Kostant
\cite{Ko}). Thus, there exists $x\in\z_\gg(e)$ such that $\A_{f+x}$
belongs to the family $\A_{t\Ad(g)\xi}$. Since $x\in\nn$, we have
$$
\lim\limits_{t\to\infty}t^2\Ad\exp(-th)(f+x)=f,
$$
and therefore $\A_f$ is a limit point of the family
$\A_{t\Ad(g)\xi}$.
\end{proof}

The "most degenerate" subalgebra $\A_\mu\subset U(\gg)$ among the
those corresponding to regular $\mu$ is the subalgebra $\A_f$. It
is a free commutative algebra with generators $\Pi_i^{(n)}$ such
that $\gr\Pi_i^{(n)}=\partial^n_e P_i\in S(\gg)$, where
$i=1,\dots,l,\ n=0,1,\dots,d_i=\deg P_i-1$, $P_i$ are the
generators of $S(\gg)^{\gg}$ (see formula \eqref{ZS}).

The element $h$ defines the \emph{principal gradation} $
\deg_{\on{pr}}$ on $U(\gg)$ such that $$\deg_{\on{pr}} e_\alpha = -
\deg_{\on{pr}} e_\alpha = 1, \quad \alpha\in\Pi,\qquad \deg_{\on{pr}}
h = 0, \quad h\in\hh.
$$
The generators of $\A_f$ are homogeneous with respect to this
gradation, with $\deg_{\on{pr}} \Pi_i^{(n)} = - n$. Thus, the algebra
$\A_f$ is graded by the principal gradation:
$\A_f=\bigoplus_{n \geq 0} \A_e^{(-n)}$. Note that the
Poincar\'e series of $\A_f$ with respect to the principal gradation is
equal to that of the algebra $U(\nn_-)$. Thus, it is natural to expect
that irreducible highest weight $\gg$-modules are cyclic as
$\A_f$-modules (having the highest weight vector as a cyclic vector).

\subsection{Opers on the projective line} Now let us describe the
spectra of Gaudin algebras following \cite{FFTL}.

Consider the Langlands dual Lie algebra ${}^L \gg$ whose Cartan
matrix is the transpose of the Cartan matrix of $\gg$. By $^L G$ we
denote the group of inner automorphisms of $^L \gg$. We fix a Cartan
decomposition
$$
^L \gg={}^L \nn\oplus {}^L \hh\oplus {}^L \nn_-.
$$
The Cartan subalgebra ${}^L \hh$ is naturally identified with
$\hh^*$. We denote by ${}^L \Delta$, ${}^L \Delta_+$, and ${}^L \Pi$
the root system of ${}^L \gg$, the set of positive roots, and the set
of simple roots, respectively.

Set
$$
p_{-1} = \sum_{\alpha^\vee\in\Pi^\vee} e_{-\alpha^\vee}\in {}^L
\gg.
$$
Let
$$
\rho=\frac{1}{2}\sum_{\alpha\in\Delta_+}\alpha\in\hh^*={}^L
\hh.
$$
The operator $\Ad \rho$ defines the principal gradation on ${}^L
\gg$, with respect to which we have a direct sum decomposition
${}^L \bb = \bigoplus_{i\geq 0} ^L \bb_i$.  Let $p_1$ be the
unique element of degree 1 in ${}^L \nn$ such that $\{
p_{-1},2\rho,p_1 \}$ is an $\ss\ll_2$-triple (note that $p_{-1}$
has the degree $-1$). Let
$$
V_{\can} = \bigoplus_{i=1}^\ell V_{\can,i}
$$
be the space of $\Ad p_1$-invariants in ${}^L \nn$, decomposed
according to the principal gradation. Here $V_{\can,i}$ has degree
$d_i$, the $i$th exponent of $^L \gg$ (and of $\g$). In particular,
$V_{\can,1}$ is spanned by $p_1$. Now choose a linear generator $p_j$
of $V_{\can,j}$.

Consider the Kostant slice in ${}^L \gg$, $$ ^L
\gg_{\can} = \left\{ p_{-1} + \sum_{j=1}^\ell y_j p_j, \quad y \in
\cc \right\}.
$$
By \cite{Ko}, the adjoint orbit of any regular element in the Lie
algebra $ ^L \gg$ contains a unique element which belongs to $^L
\gg_{\can}$. Thus, we have an isomorphism ${}^L \gg_{\can}
\tilde\to ^L \gg/^L G = \hh/W = \gg/G$.

In \cite{FF}, the center $Z(\ghat)$ of the completed enveloping
algebra $\wt{U}_{\ka_c}(\ghat)$ at the critical level is identified
with the algebra of polynomial functions on the space $\Op_{^L
G}(D^\times)$ of $^L G$-opers on the disc $D = \Spec \cc\ppart$.

The notion of opers was introduced in \cite{BD}. We refer the reader
to \cite{FFTL} for details. Here we will only say that for
$U = \Spec R$ and some coordinate $t$ on $U$, the space $\Op_{^L
G}(U)$ of $^L G$-opers is the quotient of the space of $^L
G$-connections of the form
$$
d + (p_{-1} + {\mathbf v}(t))dt, \qquad {\mathbf v}(t) \in
{}^L\bb(R)
$$
by the action of the group $^L N(R)$.

Following \cite{FFTL}, we denote by $\Op_{^L
G}(\PP^1)_{(z_i);\pi(-\mu)}$ the space of $^L G$-opers on
$\PP^1\backslash\{z_1,\dots z_N,\infty\}$ with regular
singularities at the points $z_i, i=1,\ldots,N$, and with
irredular singularity of order $2$ at the point $\infty$ with the
$2$-residue $\pi(-\mu) \in\ ^L\g/^L G = \g^*/G$, where $\pi: \g^*
\to \g^*/G$ is the porjection. Each oper from this space may be
uniquely represented in the following form:
$$
d + \left( p_{-1} + \sum_{j=1}^\ell \overline\mu_j p_j + \sum_{i=1}^N
\sum_{j=1}^\ell \sum_{n=0}^{d_j} u_{j,n}^{(i)} (t-z_i)^{-n-1} p_j
\right)dt,
$$
where
$$
p_{-1}+\sum_{j=1}^\ell \overline{\mu}_j p_j
$$
is the unique element of the $^L G$-orbit of $\mu \in \g$
contained in $^L \gg_{\can}$. Thus, $\Op_{^L
G}(\PP^1)_{(z_i);\pi(-\mu)}$ is an affine space of dimension
$\frac{1}{2}(\dim\gg+\rk\gg)N$.

\subsection{Gaudin algebra and opers}

By Theorem~5.7,(4) of \cite{FFTL}, the algebra
$\A_\mu(z_1,\dots,z_N)$ is isomorphic to a quotient of the algebra
of polynomial functions on the space $\Op_{^L
G}(\PP^1)_{(z_i);\pi(-\mu)}$, which is a polynomial algebra in
$\frac{1}{2}(\dim\gg+\rk\gg)N$ generators. On the other hand,
Proposition~\ref{predel_z} and Lemma~\ref{free pol} imply that the
algebra $\A_\mu(z_1,\dots,z_N)$ is a free polynomial algebra in
the same number of generators (see also \cite{Ryb1}, Theorem~2 and
Corollary~4). Indeed, due to the same reason as in
Lemma~\ref{predel_mu}, it suffices to prove this for the principal
nilpotent $\mu=f$. The principal gradation on $U(\gg)^{\otimes N}$
determines a filtration on $\A_f(z_1,\dots,z_N)$ such that
$\gr\A_f(z_1,\dots,z_N)=\A_f^{\otimes N}$ (since
$\lim\limits_{s\to\infty}\A_{s\mu}(z_1,\dots,z_N)=
\A_{\mu}^{(1)}\otimes\dots\otimes\A_{\mu}^{(N)}\subset
U(\gg)^{\otimes N}$, by Proposition~\ref{predel_z}). Thus $\tr
\deg \A_f(z_1,\dots,z_N)=\tr \deg \A_f^{\otimes
N}=\frac{1}{2}(\dim\gg+\rk\gg)N$.

Thus, we obtain the following assertion which was conjectured in
\cite{FFTL}, Conjecture~3.

\begin{Proposition}\label{multi-point-isomorphism} There is an isomorphism
$$
\A_\mu(z_1,\dots,z_N)\simeq\cc[\Op_{^L G}(\PP^1)_{(z_i);\pi(-\mu)}].
$$
\end{Proposition}

In particular, the algebra $\A_\mu$ is identified with the algebra of
polynomial functions on the space $\Op_{^L G}(\PP^1)_{\pi(-\mu)}$ of
$^L G$-opers on $\PP^1$ with regular singularity at the point $0$ and
with singularity of order $2$ at $\infty$, with $2$-residue
$\pi(-\mu)$, where $\pi(-\mu)$ is the image of $-\mu\in {}^L \gg$ in $
{}^L \gg/{}^L G = \gg/G$ (see \cite{FFTL}, Theorem 5.8). This space
has the following realization (see \cite{FFTL}, Section~5.4).

On the punctured disc $D_\infty^\times$ at $\infty$ (with the
coordinate $s=t^{-1}$) each element of $\Op_{^L G}(\PP^1)_{\pi(-\mu)}$
may be uniquely represented by a connection of the form
$$
d - \left(p_{-1} - \sum_{j=1}^\ell (s^{-2d_j}\overline{\mu}_j +
s^{-2d_j-1} u_j(s)) p_j\right)ds, \qquad u_j(s) = \sum_{n=0}^{d_j}
u_{j,n} s^n.
$$

On the punctured disc $D_0^\times$ at $0$, each element of $\Op_{^L
G}(\PP^1)_{\pi(-\mu)}$ may be represented uniquely by a connection of
the form
$$
d + \left( p_{-1} + \sum_{j=1}^\ell \sum_{n=0}^{d_j}
(\overline{\mu}_j + u_{j,n} t^{-n-1}) p_j \right)dt.
$$
The $1$-residue at $0$ of this oper is equal to
$$
p_{-1} + \sum_{j=1}^\ell (u_{j,d_j} + \frac{1}{4} \delta_{j,1}) p_j
\in {}^L \gg_{\can} \simeq {}^L \gg/{}^L G = \hh/W= \gg/G.
$$
In particular, for the algebra $\Spec\A_f$, we have the following
space of opers:
\begin{align*}
\Op_{^L G}(\PP^1)_{0} &= \left\{ d - \left( \left.p_{-1} -
\sum_{j=1}^\ell s^{-2d_j-1} u_j(s) p_j) ds,\ \right|\ u_j(s) =
\sum_{n=0}^{d_j} u_{j,n} s^n \right) ds \right\} \\ &= \left \{ d +
\left( p_{-1} + \sum_{j=1}^\ell \sum_{n=0}^{d_j} u_{j,n} t^{-n-1} p_j
\right)dt \right\}
\end{align*}
(here, as before, we omit $z_1$, which is set to $0$).

Next, it is proved in \cite{FFTL}, Theorem, 5.7, that for any
collection of $\gg$-modules $M_i$ with highest weights $\l_i$,
$i=1,\dots,N$, the natural homomorphism
$$
\A_\mu(z_1,\dots,z_N)\to\End(M_1\otimes\dots\otimes M_N)
$$
factors through the algebra of functions on the subspace $\Op_{^L
G}(\PP^1)_{(z_i);\pi(-\mu)}^{(\l_i)}\subset\Op_{^L
G}(\PP^1)_{(z_i);\pi(-\mu)}$ which consists of the opers with the
$1$-residue $\pi(-\l_i-\rho)$ at $z_i$. Moreover, for integral
dominant $\l_i$ the action of $\A_\mu(z_1,\dots,z_N)$ on the tensor
product of the finite-dimensional modules $V_{\l_i}$ factors through
the algebra of functions on \emph{monodromy-free} opers from $\Op_{^L
G}(\PP^1)_{(z_i);\pi(-\mu)}^{(\l_i)}$.

For every integral dominant weight $\l$, the set of monodromy-free
opers on the punctured disc $D_z^\times$ at $z$ with the regular
singularity with the residue $-\l-\rho$ at $z$ is defined by finitely
many polynomial relations. Namely, each element of the space
$\Op(D_z)^\l$ of opers on $D_z^\times$ with regular singularity and
residue $\pi(-\l-\rho)$ may be uniquely represented as
$$
d + \left(p_{-1} + \sum_{j=1}^\ell \sum_{n=-\infty}^{d_j} u_{j,n}
t^{-n-1} p_j \right) dt,
$$
with
$$
p_{-1} + \sum_{j=1}^\ell u_{j,d_j} \in\Ad({}^L G)(-\l-\rho).
$$
One can bring this connection to the form $$ d + \left(
\sum_{\alpha^\vee\in\Pi^\vee} t^{\langle\alpha^\vee,\l \rangle}
e_{-\alpha^\vee} + {\mathbf v}(t)\right)dt, $$ where ${\mathbf v}(t)
\in t^{-1}{}^L\nn+{}^L\bb[[t]]$. This oper is monodromy-free if and
only if ${\mathbf v}(t) \in {}^L\bb[[t]]$. Thus, the set of
monodromy-free opers is defined by $\dim {}^L\nn$ polynomial relations
$P_\alpha(u_{j,n})$ enumerated by positive roots
$\alpha\in\Delta_+$. These polynomial relations have the degrees
$(\alpha^{\vee},\lambda+\rho)$ with respect to the $\Z$-grading
defined by the formula $\deg u_{j,n}=-n+j$ (see \cite{FG},
Section~2.9, for details).

Introduce a $\Z$-grading on the algebra
$$\A_\mu(z_1,\dots,z_N)^{(\l_i)}=\cc \left[ \Op_{^L
G}(\PP^1)_{(z_i);\pi(-\mu)}^{(\l_i)} \right] =
\cc[u_{j,n}^{(i)}]_{j=1,\ldots,\ell; \ n=0,\dots,d_j ; \
i=1,\dots,N}$$ by the formula $\deg u_{j,n}^{(i)}=-n+j$. By comparing
the degrees of the generators of the polynomial algebras $\A_f^\l$ and
$\cc[\Op_{^L G}(\PP^1)_{0}^\l]$, we obtain the following (note that we
abbreviate the notation $\Op_{^L G}(\PP^1)_{(0);\pi(-\mu)}^\l$ to
$\Op_{^L G}(\PP^1)_{0}^\l$ when $\mu=f$):

\begin{Lemma} For $\A_f^\l=\cc[\Op_{^L G}(\PP^1)_{0}^\l]$ this grading
coincides with the principal grading on $\A_f^\l$.
\end{Lemma}


For any collection of integral dominant weights $\l_1,\dots,\l_N$
attached to the points $z_1,\dots,z_N$, we denote by
$P_{z_k;\alpha}^{(z_1,\dots,z_N);\mu;(\l_i)}$ the
polynomial in $u_{j,n}^{(i)}$ expressing the "no-monodromy" condition
at $z_k$ corresponding to the root $\alpha \in \Delta_+$.

\begin{Lemma}\label{graded}
$P_{z_k;\alpha}^{(z_1,\dots,z_N);\mu;(\l_i)}$ is an
inhomogeneous polynomial of highest degree $\deg
P_{z_k;\alpha}^{(z_1,\dots,z_N);\mu;(\l_i)}
=(\alpha^{\vee},\lambda+\rho)$. The leading term of
$P_{z_k;\alpha}^{(z_1,\dots,z_N);\mu;(\l_i)}$ is equal
to $P_{z_k;\alpha}^{(z_k);f;(\l_k)}$.
\end{Lemma}

\begin{proof} Written in terms of a local coordinate $t$ at $z_k$, an
oper from the set $\Op_{^L
G}(\PP^1)_{(z_i);\pi(-\mu)}^{(\l_i)}$ has the form
$$d + \left( p_{-1} + \sum_{j=1}^\ell \sum_{n=-\infty}^{d_j} v_{j,n}
t^{-n-1} p_j \right) dt, \quad v_{j,n}=u_{j,n}^{(k)}+\text{lower
terms}.$$ More precisely, $v_{j,n}=u_{j,n}^{(k)}$ for $n\ge0$, and for
$n<0$, the polynomial $v_{j,n}$ is a linear combination of scalars
and $u_{j,m}^{(i)},\ m\ge0$. Therefore, for $n<0$, we have $\deg
v_{j,n}<-n+j$. Hence the assertion.
\end{proof}

According to Lemma~\ref{graded}, in order to show that the
no-monodromy conditions define a finite set of opers, it suffices to
prove this for the space $\Op_{^L G}(\PP^1)_{0}^\l$. This will be done
in the next section.

\section{Main results}\label{sect-results}

\subsection{Formulation of the Main Theorem} In
section~\ref{sect-Sibuya}, we shall prove the following result.

\begin{Theorem}\label{regular_sequence} The set
of monodromy-free opers from $\Op_{^L G}(\PP^1)_{0}^\l$ is
$0$-dimensional (equivalently, the trivial monodromy conditions $
P_{0;\alpha}^{(0);f;(\l)}(u_{j,n})$ form a regular sequence).
\end{Theorem}

\subsection{Corollaries} First, let us discuss some corollaries of
this Theorem.

\begin{Corollary}\label{multi-point-generalization} The set of
  monodromy-free opers from $\Op_{^L
G}(\PP^1)_{(z_i);\pi(-\mu)}^{(\l_i)}$ is $0$-dimensional.
\end{Corollary}

\begin{proof}
This follows directly from Lemma~\ref{graded}. Indeed, the set of
monodromy-free opers from $\Op_{^L
G}(\PP^1)_{(z_i);\pi(-\mu)}^{(\l_i)}$ is the set of
common zeros of the polynomials
$P_{z_k;\alpha}^{(z_1,\dots,z_N);\mu;(\l_i)}$. By
Lemma~\ref{graded}, the leading terms of these polynomials are
$P_{z_k;\alpha}^{(z_k);f;(\l_k)}$. By
Theorem~\ref{regular_sequence}, the set of common zeros of
$P_{z_k;\alpha}^{(z_k);f;(\l_k)}$ is $0$-dimensional. Hence he set
of common zeros of the polynomials
$P_{z_k;\alpha}^{(z_1,\dots,z_N);\mu;(\l_i)}$ is also
$0$-dimensional.
\end{proof}

\begin{Corollary}\label{monodr_cond} The $\gg$-module $V_\l$ is cyclic
as an $\A_f$-module. The annihilator of $V_\l$ in $\A_{f}$ is the
ideal $I_\l\subset\A_f=\cc[\Op_{^L G}(\PP^1)_{0}^\l]$ generated by
the no-monodromy conditions on opers from $\Op_{^L G}(\PP^1)_{0}^\l$.
\end{Corollary}

\begin{proof} Note first that this assertion agrees with the
$q$-analog of the Weyl dimension formula. Namely, the Poincar\'e
series of any irreducible finite-dimensional $\gg$-module $V_\l$ with
respect to the principal grading is
$$
\chi_\l(q)=\prod\limits_{\alpha>0}\frac{1-q^{(\alpha^{\vee},
    \lambda+\rho)}}{1-q^{(\alpha^{\vee},\rho)}}.
$$
We note that the non-central generators of $\A_{f}$ have the
degrees $(\alpha^{\vee},\rho)$ with respect to the principal
grading, and the no-monodromy relations have the degrees
$(\alpha^{\vee},\lambda+\rho)$. Theorem~\ref{regular_sequence}
implies that the algebra $\A_f/I_\l$ is Gorenstein (i.e. its socle
is one-dimensional), and has the same Poincar\'e series with
respect to the principal grading as $V_\l$. Therefore the module
$V_\l$ is free as an $\A_f/I_\l$-module if and only if each
nonzero element of the socle of $\A_f/I_\l$ sends the highest
vector to some nonzero vector (which is proportional to the lowest
weight vector). Thus it remains to show that there exists an
element $a\in\A_f$ such that $av_\l=v_{w_0\l}$, (where $w_0\in W$
is the longest element of the Weyl group and $v_{w_0\l}$ is the
lowest weight vector).

Let $e,\ h,\ f$ be the principal $\sw_2$-triple containing $f$. The
module $V_\l$ decomposes into the direct sum of irreducible
$\sw_2$-modules with respect to this $\sw_2$-triple. Let $U$ be the
irreducible $\sw_2$-submodule containing $v_\l$ (this is an
$\sw_2$-submodule with the highest weight $\langle h,\l\rangle$).
Since $v_{w_0\l}$ is the unique vector of the weight $-\langle
h,\l\rangle$ with respect to the principal $\sw_2$, $U$ contains
$v_{w_0\l}$ as well. This means that we can take $a=f^{\dim
U-1}\in\A_f$.
\end{proof}

\begin{Corollary}\label{regular-elements} For any regular
$\mu\in\gg^*=\gg$, the subalgebra $\A_\mu$ has a cyclic vector in
$V_\l$. The annihilator of $V_\l$ in $\A_\mu$ is generated by the
no-monodromy conditions. In particular, the joint eigenvalues of
$\A_\mu$ in $V_\l$ (without multiplicities) are in one-to-one
correspondence with monodromy-free opers from $\Op_{^L
G}(\PP^1)_{\pi(-\mu)}^{\l}$.
\end{Corollary}

\begin{proof}  Consider the family of commutative subalgebras
$\A_{t\Ad(g)\mu}\subset U(\gg)$. By Lemma~\ref{predel_mu}, the
subalgebra $\A_f\subset U(\gg)$ is contained in the closure of this
family. Note that the condition that annihilator of $V_\l$ in
$\A_{t\Ad(g)\mu}$ is generated by the no-monodromy conditions, as well
as the existence of a cyclic vector is an open condition on
$\A_{t\Ad(g)\mu}\subset U(\gg)$. By Corollary~\ref{monodr_cond} the
subalgebra $\A_f$ satisfies both of these conditions, therefore the
conditions are satisfied for some $\A_{t\Ad(g)\mu}$. Since the
subalgebras $\A_\mu$ and $\A_{t\Ad(g)\mu}$ are conjugate, the
assertion is true for $\A_\mu$ as well. This implies that the image of
$\A_\mu$ in $\End(V_\l)$ is isomorphic to $\cc[\Op_{^L
G}(\PP^1)_{\pi(-\mu)}^{\l}]$. Hence the joint eigenvalues of $\A_\mu$
in $V_\l$ (without multiplicities) are in one-to-one correspondence
with points of $\Op_{^L G}(\PP^1)_{\pi(-\mu)}^{\l}$. (Note that this
statement was conjectured in \cite{FFTL}, Conjecture 2.)
\end{proof}

\begin{Corollary}\label{simple-spectrum} For generic $\mu\in\gg^*$ and
  any dominant integral $\l$, the quantum shift of argument subalgebra
$\A_\mu \subset U(\gg)$ is diagonalizable and has simple spectrum on
the $\g$-module $V_\l$. Moreover, its joint eigenvalues (and hence
eigenvectors, up to a scalar) are in one-to-one correspondence with
monodromy-free opers from $\Op_{^L G}(\PP^1)_{\pi(-\mu)}^{\l}$.
\end{Corollary}

\begin{proof}
First of all, by Lemma~\ref{hermitian}, the algebra $\A_\mu$ with real
$\mu$ acts by Hermitian operators on $V_\la$, and hence is
diagonalizable. By Corollary~\ref{regular-elements}, for regular $\mu$
it has a cyclic vector. The two properties may only be realized if
${\mc A}_\mu$ has simple spectrum. Hence $\A_\mu$ has simple spectrum
for regular real $\mu$. Since the simple spectrum condition is open,
$\A_\mu$ has simple spectrum for generic $\mu \in \gg^*$. By
Corollary~\ref{regular-elements}, the joint eigenvalues of $\A_\mu$ in
$V_\l$ are in one-to-one correspondence with monodromy-free opers from
$\Op_{^L G}(\PP^1)_{\pi(-\mu)}^{\l}$.
\end{proof}

\begin{Corollary}\label{regular-elements-z}
For any $N$-tuple of pairwise distinct complex numbers
$z_1,\dots,z_N\in\cc$ and any regular $\mu\in\gg^*=\gg$, the
subalgebra $\A_\mu(z_1,\dots,z_N)\subset U(\gg)^{\otimes N}$ has a
cyclic vector in $V_{(\l_i)}=V_{\l_1}\otimes\dots\otimes V_{\l_n}$.
The annihilator of $V_{(\l_i)}$ in $\A_\mu(z_1,\dots,z_N)$ is
generated by the no-monodromy conditions.
\end{Corollary}

\begin{proof} Due to the same reason as in
Corollary~\ref{regular-elements}, it suffices to prove this for
the principal nilpotent $\mu=f$. Let $e,h,f$ be the principal
$sl_2$-triple containing $f$, and let $s\in\cc$. Then $$\exp(\ad
sh)(\A_f(z_1,\dots,z_N))=\A_{\exp(-2s)f}(z_1,\dots,z_N).$$ We have
$$\lim\limits_{s\to-\infty}\A_{\exp(-2s)f}(z_1,\dots,z_N)=
\A_{f}^{(1)}\otimes\dots\otimes\A_{f}^{(N)},$$ by
Proposition~\ref{predel_z}). By Corollary~\ref{regular-elements},
the algebra $\A_{f}^{(1)}\otimes\dots\otimes\A_{f}^{(N)}$ has a
cyclic vector in $V_{(\l_i)}=V_{\l_1}\otimes\dots\otimes V_{\l_n}$.
Hence $\A_{\exp(-2s)f}(z_1,\dots,z_N)$ has a cyclic vector for
some $s\in\cc$. Since $$\exp(\ad
sh)(\A_f(z_1,\dots,z_N))=\A_{\exp(-2s)f}(z_1,\dots,z_N),$$ the
algebra $\A_f(z_1,\dots,z_N)$ has a cyclic vector as well.

The assertion on the annihilator of $V_{(\l_i)}$ is proved by the same
reasoning as in Corollary~\ref{regular-elements} with the reference to
Lemma~\ref{graded}. (Note that this implies Conjecture~4 of
\cite{FFTL}.)
\end{proof}

\begin{Corollary}\label{simple-spectrum_z}
For generic $z_1,\dots,z_N\in\cc,\ \mu\in\gg^*$, the subalgebra
$\A_\mu(z_1,\dots,z_N)\subset U(\gg)^{\otimes N}$ has simple
spectrum in $V_{(\l_i)}=V_{\l_1}\otimes\dots\otimes V_{\l_N}$. Hence
the joint eigenvectors for higher Gaudin hamiltonians in
$V_{(\l_i)}$ are in one-to-one correspondence with monodromy-free
opers from $\Op_{^L
G}(\PP^1)_{(z_i);\pi(-\mu)}^{(\l_i)}$.
\end{Corollary}

\begin{proof}
The same reasoning as in the proof of
Corollary~\ref{simple-spectrum}.
\end{proof}

\subsection{Idea of the Proof of the Main
Theorem}\label{subsect-idea} The proof of
Theorem~\ref{regular_sequence} will be given in the next section. The
main tools we use are the following local normalization theorems for
irregular singular connections.

\begin{Fact} \emph{(Hukuhara--Turritin--Levelt theorem, see \cite{BV},
    and also \cite{IYa,Va} for the ${\mf g}{\mf l}_n$ case)} Any
connection $d+A(t)dt$ on the punctured disc $D^\times =\on{Spec}
\cc\ppart$, where $A(t)=\sum\limits_{k=-2}^\infty A_kt^k$ may be
reduced by a suitable \emph{formal} shearing gauge transformation
$H(t^{\frac{1}{N}})$ to its \emph{formal normal form}
\begin{equation}\label{Hukuhara-Turritin-Levelt} d+B(w)dw\quad
B(w)=\sum\limits_{k=2}^m B_kw^{-k}+Cw^{-1},\end{equation} where
$w^N=t$, $m\le N+1$ and all $B_k$ commute with $C$ and belong to
a fixed Cartan subalgebra. Moreover, if $A_{-2}$ is nilpotent,
then $m\le N$.
\end{Fact}

\begin{Fact} \emph{(Sibuya sectorial normalization
theorem, see \cite{Bo}, Appendix A. For the ${\mf g}{\mf l}_n$ case,
see \cite{Si}, \cite{Was}, \cite{IYa})} For any sector $S$ of opening
$\frac{\pi}{m-1}$ on the $w$-plane, the formal gauge transformation
$H(w)$ may be extended to an {\em analytic} gauge transformation
$H_S(w)$ conjugating the connection $d+A(t)dt$ to its formal normal
form.
\end{Fact}

The main idea of the proof of Theorem~\ref{regular_sequence} is as
follows. The opers we are interested in are represented by connections
of the form $d+A(s)ds$, where $A(s)=\sum\limits_{k=-2}^\infty A_ks^k$
with {\em nilpotent} $A_{-2}$ in the neighborhood of $\infty$. By
Hukuhara--Turritin--Levelt Theorem, its formal normal form
is~(\ref{Hukuhara-Turritin-Levelt}) with $m\le N$. Moreover, all $B_k$
are zero if and only if the initial connection has regular singularity
at $\infty$, i.e. $u_{j,n}=0$ for $n\ne d_j$. The solutions of
$^L\gg$-valued ordinary differential equation
\begin{equation}    \label{formal eq}
d\phi(w)+[B(w),\psi(w)]dw = 0
\end{equation}
have the form
$$
\Ad\left(\exp\left(\sum\limits_{k=2}^m
B_kw^{-k+1}\right)w^C\right)x=\Ad\left(\exp\left(\sum\limits_{k=2}^m
B_kt^{\frac{-k+1}{N}}\right)s^\frac{C}{N}\right),
$$
where $x\in {}^L\gg, m\le N$ (exponentials of linear combinations
of \emph{fractional} powers of $s$). Each sector $S$ on the
$w$-plane with the origin at $0$, which does not contain real
multiples of the eigenvalues of the operators $\Ad \exp(\frac{2\pi
i k}{m-1})B_m$ for $k=1,\dots,m-1$, distinguishes a subspace of
formal solutions of the equation \eqref{formal eq} which
\emph{decay exponentially} as $w\to0$ along each ray in this
sector $S$. Due to the Sibuya theorem, on such a sector we also
have a subspace of exponentially decaying solutions of the
equation
\begin{equation}    \label{formal eq1}
d\phi(s)+[A(s),\phi(s)]ds = 0.
\end{equation}

We have assumed that the connection $d+A(s)ds$ has trivial monodromy
representation. Hence we may pick a global solution $\phi$ to the
$^L\gg$-valued ordinary differential equation \eqref{formal eq1} which
decays exponentially as $w\to0$ along each ray in some sector $S_0$ on
the $w$-plane. Since solution $\phi$ is a single-valued function of
$s=w^N$, we find that $\phi$ also decays on the sector $\exp(2\pi i
k/N)S_0 , k=1,\ldots,N-1$. Consider the sector $S_j = \exp(2\pi i
/N)S_0$. We will show in section~\ref{sect-Sibuya}, that there exists
a sector $S_r$ (actually $S_r=\exp(\frac{\pi i}{m-1}) S_0$) such that
\begin{enumerate} \item there is a formal solution of \eqref{formal
eq1} which decays simultaneously on $S_0$ and $S_j$; \item there
is no formal solution of \eqref{formal eq1} which decays
simultaneously on $S_r$ and $S_0$; \item there is a sector $S$ on
the $w$-plane of the opening $\frac{\pi}{m-1}$ which has nonzero
intersection with $S_j$ and $S_r$; \item there is a sector $S'$ on
the $w$-plane of the opening $\frac{\pi}{m-1}$ which has nonzero
intersection with $S_0$ and $S_r$.
\end{enumerate}
Applying Sibuya's sectorial normalization theorem to the sectors
$S$ and $S'$, we obtain that there exists a non-zero formal
solution of \eqref{formal eq1} exponentially decaying on both
$S_0$ and $S_r$, which contradicts the second condition.

Thus, for any monodromy-free oper from $\Op_{^L G}(\PP^1)_0^\l$,
we have $u_{j,n}=0$ for $n\ne j$ and $u_{j,d_j}$ are fixed by the
residue at $0$. Hence the space of monodromy-free opers from
$\Op_{^L G}(\PP^1)_0^\l$ is $0$-dimensional. This completes the
proof.

\section{Proof of the Main Theorem}\label{sect-Sibuya}

In this section we prove Theorem~\ref{regular_sequence}.

\subsection{The ${\mf s}{\mf l}_2$ case}

Let $d+A(t)dt$ be a ${\mf s}{\mf l}_2$-connection, such that
\begin{equation}\label{gl2-irreg-nilp}
  A(t)=e_{21}+(\l(\l+1)t^{-2}+a^2t^{-1})e_{12},\end{equation}
where $\l$ is a fixed integral weight, $a\in\cc$.

\begin{Proposition}
The connection~(\ref{gl2-irreg-nilp}) has trivial monodromy
representation if and only if $a=0$.
\end{Proposition}

\begin{proof}
Let $s=t^{-1}$ be the local coordinate at the infinity. We rewrite
our connection in $s$ as follows.
$$A(s)=s^{-2}e_{21}+(\l(\l+1)+a^2s^{-1})e_{12}.$$ Changing the
variable as $s=w^2$, we obtain the following connection:\footnote{Note
that if we were considering an oper in $\on{Op}_{^L G}(\pone)_\mu$,
where $\mu$ is regular semi-simple, then we would be able to bring it
to a normal form without extracting the square root of $s$. For this
reason the argument given below would not work in this case.}
$$d+A(w)dw,\ A(w)=2(w^{-3}f_{21}+(\l(\l+1)w+a^2w^{-1}))e_{12}.$$ The
latter is conjugate to
\begin{equation}\label{gl2-irreg}d+((2ae_{11}-2ae_{22})w^{-2}+B(w))dw,
\end{equation}
where $B(w)=O(w^{-1})$ as $w\to0$. The connection (\ref{gl2-irreg}) is
\emph{formally} conjugate to
\begin{equation}\label{gl2-irreg-formal}d+((2ae_{11}-2ae_{22})w^{-2}
+(b_{11}e_{11}+b_{22}e_{22})w^{-1})dw.\end{equation}
Moreover, by the Sibuya sectorial normalization theorem, for any
sector $S=\{w|\Arg w\in(\alpha,\alpha+\pi)\}$, the formal gauge
transformation $H(w)$ can be extended to an analytical gauge
transformation $H_S(w)$ conjugating (\ref{gl2-irreg}) to
(\ref{gl2-irreg-formal}).

Consider two sectors, $S_0=\{w|\Arg w\in(\Arg a-\frac{\pi}{2},\Arg
a+\frac{\pi}{2})\}$ and $S_1=\{w|\Arg w\in(\Arg
a+\frac{\pi}{2},\Arg a+\frac{3\pi}{2})\}$. We have the following
basis of solutions to the linear ordinary differential equation
(\ref{gl2-irreg-formal})
$$\psi_0(w)=(0,\exp(-2aw^{-1})w^{b_{22}}),\quad \psi_1(w)=
(\exp(2aw^{-1})w^{b_{11}},0).
$$
Note that the solution $\psi_0$ decays exponentially as $w\to0$
along each ray in the sector $S_0$ and blows up exponentially as
$w\to0$ along each ray in the sector $S_1$. Respectively, $\psi_1$
decays on $S_1$ and blows up on $S_0$.

Assume that the connection (\ref{gl2-irreg-nilp}) has trivial
monodromy. Let $\phi$ be the global solution to the equation
(\ref{gl2-irreg}) such that $\phi|_{S_0}=H_{S_0}\psi_0$. Since
$\psi_0$ decays exponentially on $S_0$ and the gauge transformation
$H_{S_0}(w)$ is bounded in some neighborhood of $0$, the solution
$\phi$ also decays on ${S_0}$. Since solution $\phi$ is a
single-valued function of $s=w^2$, we have $\phi(w)=\phi(-w)$, and
hence $\phi|_{S_1}$ decays as well.

Consider the sector $S=\{w|\Arg w\in(\Arg a,\Arg a+\pi)\}$. On
$S$, we have the following basis of solutions to the
equation~(\ref{gl2-irreg}):
$$\phi_0=H_S\psi_0,\quad \phi_1=H_S\psi_1.
$$
We have $\phi=k_0\phi_0+k_1\phi_1$ for some $k_0,k_1\in\cc$. Since
$\phi$ decays on $S\cap S_0$ while $\psi_1$ blows up on $S\cap
S_0$, we have $k_1=0$. Since $\phi$ decays on $S\cap S_1$ while
$\psi_0$ blows up on $S\cap S_1$, we have $k_0=0$. Hence $\phi=0$
and we have a contradiction.
\end{proof}

\subsection{General case} Let $d+A(t)dt$ be a $\gg$-connection, such
that $A(t)=p_{-1}+\sum\limits_{j=1}^\ell
c_jt^{-d_j-1}p_j+\sum\limits_{n=0}^{d_j-1}u_{j,n}t^{-n-1}p_j$.

\begin{Proposition}
The connection $d+A(t)dt$ has trivial monodromy representation if
and only if $u_{j,n}=0$ for all $j,n$.
\end{Proposition}

\begin{proof}
Suppose that the connection~$d+A(t)dt$ has no monodromy while some
of the $u_{j,n}$ are nonzero.

Let $s=t^{-1}$ be the local coordinate at the infinity. We rewrite
our connection in $s$ as follows.
\begin{equation}\label{irreg-conn-general-1}d+A(s)ds\quad
A(s)=s^{-2}p_{-1}+\sum\limits_{j=1}^\ell
c_js^{d_j-1}p_j+\sum\limits_{n=0}^{d_jj-1}u_{j,n}s^{n-1}p_j.
\end{equation}

Let us change the variable as $w^N=s$, where $N=\prod\limits_{j\in
E} j$. The connection rewrites as
\begin{equation}\label{irreg-conn-general}d+A(w)dw,\quad
A(w)=w^{-N-1}p_{-1}+\sum\limits_{j\in
E}c_jw^{jN-1}p_j+\sum\limits_{n=0}^{j-1}u_{j,n}w^{nN-1}p_j.
\end{equation}

Due to the Hukuhara--Turritin--Levelt theorem, the connection
(\ref{irreg-conn-general}) may be reduced by a suitable
\emph{formal} gauge transformation $H(w)$ to its \emph{formal
normal form}
\begin{equation}\label{irreg-conn-general-formal}d+B(w)dw\quad
B(w)=\sum\limits_{k=2}^m B_kw^{-k}+Cw^{-1},\end{equation} where $m\le
N$ and all $B_k$ commute with $C$ and belong to the fixed Cartan
subalgebra $\hh \subset \gg$. The coefficients $B_k$ are zero (i.e.,
the connection has regular singularity) if and only if all $u_{j,n}$
are zero.

Let $\Delta_+$ be the set of positive roots with respect to a Borel
subalgebra containing $\hh$. To any root $\alpha\in \Delta_+$ such
that $\alpha(B_m)\ne0$, we assign a collection of $2m-2$
\emph{separation rays} (also known as \emph{Stokes directions})
defined by the condition
\begin{equation}\label{stokes-dir} \re \frac{\alpha(B_m)}{w^{m-1}}=0.
\end{equation}

\begin{Remark} The set of Stokes directions is invariant under the
transformations $w\mapsto\exp(\frac{\pi i}{m-1})w$ (clear from the
condition~(\ref{stokes-dir})) and $w\mapsto\exp(\frac{2\pi
i}{N})w$ (since the connection~(\ref{irreg-conn-general-1})
depends on $s=w^N$, not $w$).
\end{Remark}

{\bf Generic case.} Assume that $B_m$ is \emph{generic},
i.e., $B_m$ is a regular element of $\hh$. Then all separation rays
are distinct. The proof under this assumption is slightly simpler, so
we shall give the detailed proof in this case first, and then explain
what should be changed for an arbitrary $B_m$.

In the generic case we have precisely $M=(2m-2)r$ separation rays,
where $r=|\Delta_+|$. We label these separation rays
$d_1,d_2,\dots d_M=d_0$ going counter-clockwise, and choose the
initial sector $S_0$ to be between $d_0$ and $d_1$.

Consider the $\gg$-valued ordinary differential equation
\begin{equation}\label{irreg-equation}d\phi(w)+[A(w),\phi(w)]dw=0,
\end{equation}
and its formal normal form
\begin{equation}\label{irreg-equation-formal}d\phi(w)+[B(w),\phi(w)]
dw=0.\end{equation}
We have the following basis of solutions of equation
(\ref{irreg-equation-formal}):
\begin{equation}\label{solution-basis-formal}\psi_\alpha(w)=\Ad
\exp\left(\sum\limits_{k=2}^m \frac{1}{k-1}
B_kw^{-k+1}\right)w^{-C}e_\alpha, \quad \psi_{h_j}(w)=\Ad w^{-C}h_j,
\end{equation}
where $\alpha\in {}\Delta$, and $h_j, j=1,\ldots,\ell$, is a basis of
$\hh$. On each sector $S_i$ between the separation rays $d_i$ and
$d_{i+1}$, the formal solutions behave as follows:
\begin{enumerate}\item If $\re \frac{\alpha(B_m)}{w^{m-1}}<0$ on
$S_i$, then $\psi_\alpha(w)$ decays exponentially as $w\to0$ along
each ray in $S_i$, \item If $\re \frac{\alpha(B_m)}{w^{m-1}}>0$ on
$S_i$, then $\psi_\alpha(w)$ blows up exponentially as $w\to0$ along
each ray in $S_i$, \item $\psi_{h_j}(w)$ has polynomial growth/decay
as $w\to0$ along each ray in $S_i$.
\end{enumerate}
We note that if $|i-j|<r$ then there is a root $\alpha\in
{}\Delta$ such that $\re \frac{\alpha(B_m)}{w^{m-1}}<0$ on both
$S_i$ and $S_j$. This means, that, for $|i-j|<r$, there is a
formal solution $\psi^{(ij)}=\psi_\alpha$ which decays on both
sectors $S_i$ and $S_j$.

Note that the formal solution are periodic with respect to
$w\mapsto\exp(\frac{2\pi i}{m-1})w$. We shall consider the
behavior of solutions in the "fundamental" sector $\widehat{S_0}$
between $d_M=d_0$ and $d_{2r}$. For each sector
$S_i\subset\widehat{S_0}$ except $S_0$, there is a formal solution
which decays simultaneously on $S_i$ and $S_r$. Moreover, for each
sector $S_i\subset\widehat{S_0}$ there is a sector of opening
$\frac{\pi}{m-1}$ (a half-period) which has nonzero intersection
with $S_i$ and $S_r$.

On the other hand, the global solutions of the
equation~(\ref{irreg-equation}) are periodic with respect to
$w\mapsto\exp(\frac{2\pi i}{N})w$. Since the set of separation
rays is invariant under this transformation, this transformation
permutes the sectors $S_i$. Hence the global solutions have the
same asymptotic properties on the sectors $S_0$ and
$S_j=\exp(\frac{2\pi i}{N})S_0$. Since $m-1<N$, the both sectors
$S_0$ and $S_j=\exp(\frac{2\pi i}{N})S_0$ belong to the
fundamental sector $\widehat{S_0}$.

Consider the sector $S_j=\exp(\frac{2\pi i}{N})S_0$. Since
$m-1<N$, we have $0<j<2r$, and therefore there is a formal
solution $\psi^{(jr)}$ which decays on $S_j$ and $S_r$. Moreover,
since the opening of the sector between $d_j$ and $d_r$ is less
than $\frac{\pi}{m-1}$, there is a sector $S_{jr}$ of opening
$\frac{\pi}{m-1}$ which contains $S_j$ and $S_r$. By Sibuya's
theorem, the formal gauge transformation $H(w)$ can be extended to
an analytical gauge transformation $H_{S_{jr}}(w)$ conjugating
(\ref{irreg-conn-general}) to (\ref{irreg-conn-general-formal}) on
$S_{jr}$. Thus, there is a solution
$\phi^{(jr)}:=H_{S_{jr}}\psi^{(jr)}$ of the
equation~(\ref{irreg-equation}), which decays on both sectors
$S_j$ and $S_r$. Due to the no-monodromy condition, there exists a
global solution $\phi$ such that $\phi|_{S_{jr}}=\phi^{(jr)}$.
Moreover, the no-monodromy condition means that each solution of
the equation~(\ref{irreg-equation}) is a single-valued function of
$s=w^N$, hence we have $\phi(w)=\phi(\exp(\frac{2\pi i}{N})w)$.
This means that $\phi$ decays exponentially on
$S_0=\exp(\frac{-2\pi i}{N})S_j$.

According to the Sibuya theorem, since the opening of the sector
between $d_M=d_0$ and $d_r$ is $\frac{\pi}{m-1}$, the
equation~(\ref{irreg-equation}) can be conjugated to its formal
normal form in some sector $S$ intersecting $S_0$ as well as
$S_r$. On $S$, we have the following basis of solutions to the
equation~(\ref{irreg-equation})
$$H_S\psi_\alpha,\ H_S\psi_{h_j},\quad\alpha\in {}\Delta,
\ j=1,\dots,\ell.
$$

We have
$$
\phi=\sum_{\alpha\in {}\Delta}k_\alpha
H_S\psi_\alpha+\sum_{j=1}^{\ell}k_{h_j}H_S\psi_{h_j}.
$$
For any $h_j$, the solution $H_S\psi_{h_j}$ does not decay
exponentially on $S$, hence $k_{h_j}=0$ for all
$j=1,\dots,\ell$. Next, for any $\alpha\in {}\Delta$ we have the
following alternative: {\bf either} $\re
\frac{\alpha(B_m)}{w^{m-1}}>0$ on $S_0$, and hence $H_S\psi_\alpha$
blows up on $S\cap S_0$, {\bf or} $\re \frac{\alpha(B_m)}{w^{m-1}}>0$
on $S_r$, and hence $H_S\psi_\alpha$ blows up on $S\cap S_r$. This
means that $k_\alpha=0$ for all $\alpha\in {}\Delta$. Hence $\phi=0$,
and we have a contradiction.

\bigskip

{\bf General case.} For non-generic $B_m$, the asymptotic behavior
of formal solutions is not determined only by $B_m$, but depends
also on $B_k$ with $k<m$. Thus it is difficult to figure out on
which sectors the given formal solution decays exponentially. The
crucial observation is that for our purposes it is sufficient to
watch only for the \emph{most rapid} decay of solutions (i.e.
faster than $\exp(bw^{-m+1})$ for some $b$), which is determined
by the leading term $B_m$.

We have $M=(2m-2)r$ separation rays, where
$r=\#\{\alpha|\alpha(B_m)\ne0\}$. We label the separation rays
$d_1,\dots d_M$ going in the positive sense, and choose the
initial sector $S_0$ between $d_M$ and $d_1$. Note that, for
non-generic $B_m$, some of the separation rays may coincide. We
choose the initial sector $S_0$ to be non-empty.

Consider the $ \gg$-valued ordinary differential
equation~(\ref{irreg-equation}), and its formal normal
form~(\ref{irreg-equation-formal}). We have the
basis~(\ref{solution-basis-formal}) of solutions to the equation
(\ref{irreg-equation-formal}), which behaves on each sector $S_i$
between the separation rays $d_i$ and $d_{i+1}$ as follows:
\begin{enumerate}\item If $\re \frac{\alpha(B_m)}{w^{m-1}}<0$ on
  $S_i$, then $\psi_\alpha(w)$
decays \emph{most rapidly} (i.e. faster than $\exp(bw^{-m+1})$ for
some $b$) as $w\to0$ along each ray in $S_i$, \item If $\re
\frac{\alpha(B_m)}{w^{m-1}}>0$, then $\psi_\alpha(w)$ blows up as
$w\to0$ along each ray in $S_i$,\item If $\re
\frac{\alpha(B_m)}{w^{m-1}}=0$, then $\psi_\alpha(w)$ does not
decay or decays not faster than $\exp(bw^{-m+1})$ for all
$b\in\cc$ as $w\to0$ along each ray in $S_i$, \item
$\psi_{h_j}(w)$ has polynomial growth/decay as $w\to0$ along each
ray in $S_i$.
\end{enumerate}
As in the generic case, we note that if $|i-j|<r$ then
there is a root $\alpha\in {}\Delta$ such that $\re
\frac{\alpha(B_m)}{w^{m-1}}<0$ on both $S_i$ and $S_j$, and
hence, there is a formal solution $\psi^{(ij)}=\psi_\alpha$ which
decays most rapidly on both sectors $S_i$ and $S_j$.

In the same way as in the generic case, we get the global solution
$\phi$ of the equation~(\ref{irreg-equation}) which decays most
rapidly on $S_0$ and $S_r$. We take a sector $S$ intersecting both
$S_0$ and $S_r$, and represent the solution $\phi$ on $S$ as
$$
\phi=\sum_{\alpha\in {}\Delta}k_\alpha
H_S\psi_\alpha+\sum_{j=1}^{\ell}k_{h_j}H_S\psi_{h_j}.
$$
For any $h_j$, the solution $H_S\psi_{h_j}$ does not decay
exponentially on $S$, hence $k_{h_j}$ for all $j=1,\dots,\ell$. For
any $\alpha\in {}\Delta$, we have the following possibilities:
\begin{enumerate} \item $\re \frac{\alpha(B_m)}{w^{m-1}}>0$ on $S_0$,
and hence $H_S\psi_\alpha$ blows up on $S\cap S_0$, \item $\re
\frac{\alpha(B_m)}{w^{m-1}}>0$ on $S_r$, and hence
$H_S\psi_\alpha$ blows up on $S\cap S_r$, \item $\re
\frac{\alpha(B_m)}{w^{m-1}}=0$, and hence $H_S\psi_\alpha$ does
not decay or decays not faster than $\exp(bw^{-m+1})$ for all
$b\in\cc$ on $S$.\end{enumerate} This means that $k_\alpha=0$ for
all $\alpha\in {}\Delta$. Hence $\phi=0$, and we have a
contradiction.

\end{proof}

\end{document}